\documentclass[a4paper,11pt]{amsproc}

\usepackage{amsmath,amsthm,amssymb,mathtools,mathrsfs,fullpage,hyperref}

\usepackage{xcolor}
\usepackage[capitalize,nameinlink,noabbrev,nosort]{cleveref}
\hypersetup{
    colorlinks=true,       
    linkcolor=blue,          
    citecolor=blue,        
    urlcolor=blue,           
}

\makeatletter
\@namedef{subjclassname@2020}{%
  \textup{2020} Mathematics Subject Classification}
\makeatother

\theoremstyle{plain}
\newtheorem*{theorem}{Theorem}

\newcommand{\Z}{\mathbb{Z}}

\begin{document}
\title[]{On the length over which $k$-G\"obel sequences remain integers} 

\author{Yuh Kobayashi}
\address{Department of Mathematical Sciences, Aoyama Gakuin University, 5-10-1 Fuchinobe, Chuo-ku, Sagamihara, Kanagawa, 252-5258, Japan}
\email{kobayashi@math.aoyama.ac.jp}

\author{Shin-ichiro Seki}
\address{Department of Mathematical Sciences, Aoyama Gakuin University, 5-10-1 Fuchinobe, Chuo-ku, Sagamihara, Kanagawa, 252-5258, Japan}
\email{seki@math.aoyama.ac.jp}

\subjclass[2020]{Primary 11B37; Secondary 11B50}
\keywords{G\"{o}bel sequences}
\thanks{This research was supported by JSPS KAKENHI Grant Numbers JP22K13960 (Kobayashi) and JP21K13762 (Seki).}
\begin{abstract}
We prove that the sequence $(N_k)_k$, where each $N_k$ is defined as the smallest positive integer $n$ for which the $n$th term $g_{k,n}$ of the $k$-G\"obel sequence is not an integer, is unbounded.
\end{abstract}
\maketitle
\section{Introduction}
For $k\geq 2$, the \emph{$k$-G\"obel sequence} $(g_{k,n})_n$ is defined by the initial value $g_{k,1}=2$ and the recursion $ng_{k,n}=(n-1)g_{k,n-1}+g_{k,n-1}^k$.
Let $N_k\coloneqq\inf\{n\geq 1 \mid g_{k,n}\not\in\Z\}$.
The original G\"obel's sequence (= 2-G\"obel sequence, \cite[A003504]{Sloane}) with $N_2=43$ has attracted interest as an example of \emph{the strong law of small numbers} (\cite{Guy}). (The growth of $g_{k,n}$ is very fast.
In fact, the value $g_{2,43}\approx5.4\times 10^{178485291567}$ is very large.)

The behavior of the sequence $(N_k)_k$ (\cite[A108394]{Sloane}) is not yet understood well and remains mysterious.
In \cite{MatsuhiraMatsusakaTsuchida2024}, Matsuhira, Matsusaka, and Tsuchida proved that $\min_{k\geq 2} N_k = 19$.
As mentioned in \cite[Section~3]{MatsuhiraMatsusakaTsuchida2024} and \cite[Episode~3]{KobayashiDSeki}, the following three questions are fundamental problems about the sequence $(N_k)_k$:
\begin{enumerate}
\item\label{it1} Why is $N_k$ a prime number for most values of $k$? Or rather, in what cases does $N_k$ become a composite number?
\item\label{it2} Does $N_k$ always take a finite value for any $k\geq 2$?
\item\label{it3} Is the sequence $(N_k)_k$ unbounded?
\end{enumerate}
For Problems~\eqref{it1} and \eqref{it2}, only numerical data has been obtained.
In \cite{KobayashiYSeki}, it is shown that $N_k$ is prime for 86.5\% of values up to $k\leq 10^7$, and that $N_k$ is finite for $k\leq 10^{14}$.

In this short note, however, we report that Problem~\eqref{it3} can be solved in a very elementary way within the framework of \cite{MatsuhiraMatsusakaTsuchida2024}.
Let $m\#$ denote the primorial of $m$ or, in other words, the radical of $m!$.
\begin{theorem}
Let $m$ be a positive integer.
If $k\geq 2$ satisfies $k\equiv 1 \pmod{m!/m\#}$, then $N_k > m$.
In particular, $\sup_{k\geq 2}N_k=\infty$.
\end{theorem}
\section{Proof}
Let $k\geq 2$, $r\geq 1$ be integers and $p$ a prime.
Let $\Z_{(p)}$ be the localization of $\Z$ at $(p)$ and $\nu_p(n)$ be the $p$-adic valuation of $n$.
For any positive integer $n$ with $\nu_p(n!)\leq r$, we define $g_{k,p,r}(n)\in\Z/p^{r-\nu_p(n!)}\Z \cup \{\mathsf{F}\}$ as in \cite[Definition~9]{MatsuhiraMatsusakaTsuchida2024}: For $n=1$, $g_{k,p,r}(1)=2\bmod{p^r}$.
For $n\geq 2$, when $g_{k,p,r}(n-1)=\mathsf{F}$, $g_{k,p,r}(n)=\mathsf{F}$.
When $g_{k,p,r}(n-1)=a\bmod{p^{r-\nu_p((n-1)!)}}$, if $(n-1)a+a^k\not\equiv 0\pmod{p^{\nu_p(n)}}$, then $g_{k,p,r}(n)=\mathsf{F}$, and otherwise, $g_{k,p,r}(n)=\frac{(n-1)a+a^k}{p^{\nu_p(n)}}\cdot c\bmod{p^{r-\nu_p(n!)}}$, where $c$ is an integer such that $c\cdot (n/p^{\nu_p(n)})\equiv 1\pmod{p^{r-\nu_p(n!)}}$.
Let $\varphi(n)$ denote Euler's totient function.
\begin{proof}[Proof of Theorem]
Let $m\geq 4$ and $k\geq 2$ be integers.
Assume that $k\equiv 1\pmod{m!/m\#}$.
For each prime $p\leq m$, set $r_p\coloneqq\nu_p(m!)$ and $k_p\coloneqq \varphi(p^{r_p})+1$.
It is clear that $k_p > r_p \geq 1$.
Let us temporarily suppose that for some prime $p\leq m$, we have $g_{k_p,p,r_p}(m)\neq\mathsf{F}$.
Since $\varphi(p^{r_p})$ divides $m!/m\#$, we see that $k\equiv k_p \pmod{\varphi(p^{r_p})}$, and thus, by \cite[Proposition~12]{MatsuhiraMatsusakaTsuchida2024}, it follows that $g_{k,p,r_p}(m)=g_{k_p,p,r_p}(m)\neq \mathsf{F}$.
Therefore, by \cite[Lemma~10]{MatsuhiraMatsusakaTsuchida2024}, we conclude that $g_{k,n}\in\Z_{(p)}$ for all $1\leq n\leq m$.
By this argument, together with \cite[Lemma~4]{MatsuhiraMatsusakaTsuchida2024} and the fact that $\bigcap_p\Z_{(p)}=\Z$, in order to prove $N_k > m$, it suffices to show that $g_{k_p,p,r_p}(m)\neq\mathsf{F}$ for each prime $p\leq m$.

Since $k_2>r_2$, we see that $2+2^{k_2} \equiv 2 \pmod{2^{r_2}}$.
Hence, we have $g_{k_2,2,r_2}(2)=1\bmod{2^{r_2-1}}$.
It is clear that, subsequently, $g_{k_2,2,r_2}(n)=1\bmod{2^{r_2-\nu_2(n!)}}$ for $2\leq n\leq m$.
In particular, $g_{k_2,2,r_2}(m)\neq\mathsf{F}$.

Let $p\leq m$ be an odd prime.
For any $1\leq n\leq m$, we have $g_{k_p,p,r_p}(n)=2\bmod{p^{r_p-\nu_p(n!)}}$.
In fact, if $g_{k_p,p,r_p}(n-1)=2\bmod{p^{r_p-\nu_p((n-1)!)}}$, then since $(n-1)2+2^{k_p}\equiv 2n \pmod{p^{r_p-\nu_p((n-1)!)}}$ by Euler's theorem, we have $g_{k_p,p,r_p}(n)=2\bmod{p^{r_p-\nu_p(n!)}}$.
In particular, $g_{k_p,p,r_p}(m)\neq\mathsf{F}$.
\end{proof}
\section{Remark}
By replacing the initial value in the definition of the $k$-G\"obel sequence with $g_{k,1}=l$, we define the \emph{$(k,l)$-G\"obel sequence}, which has been investigated in \cite{GimaEtal2024+,Ibstedt1990,KobayashiYSeki}.
Generalizing our previous arguments as follows, the theorem holds in the same form for $(k,l)$-G\"obel sequences as well.

Fix $l$ and a prime $p\leq m$, and in the definition of $g_{k,p,r}$, replace the initial condition with $g_{k,p,r}(1)=l\bmod{p^r}$.
We use the notation $r_p$ and $k_p$ as in the previous section.
We can easily check that by induction on $n$, for $n\leq p^{\nu_p(l)}$, $g_{k_p,p,r_p}(n)=(l/p^{\nu_p(n)})\cdot c_{n,p}\bmod{p^{r_p-\nu_p(n!)}}$ holds, while for $n\geq p^{\nu_p(l)}$, we have $g_{k_p,p,r_p}(n)=l/p^{\nu_p(l)}\bmod{p^{r_p-\nu_p(n!)}}$.
Here, $c_{n,p}\in\Z$ satisfies $c_{n,p}\cdot(n/p^{\nu_p(n)})\equiv 1 \pmod{p^{r_p-\nu_p(n!)}}$.
Therefore, the same proof works for a general $l$.

\end{document}